\documentclass[10pt]{article}

\usepackage{mathrsfs}
\usepackage{amsfonts}
\usepackage{amsfonts}
\usepackage{graphicx,latexsym,bm,amsmath,amssymb,verbatim,multicol,lscape}

\date{}

\begin{document}
\title{{\bf Counting rational points of an algebraic variety over finite fields}
\author{Shuangnian Hu$^{\rm a, b}$, Shaofang Hong$^{\rm a}$ and Xiaoer Qin$^{\rm c}$\\
$^{\rm a}$ {\it Mathematical College, Sichuan University, Chengdu 610064, P.R. China} \\
$^{\rm b}$ {\it School of Mathematics and Statistics, Nanyang Institute of Technology,}\\
{\it Nanyang 473004, P.R. China} \\
$^{\rm c}$ {\it College of Mathematics and Statistics, Yangtze Normal University,}\\
{\it Chongqing 408100, P.R. China} \\
E-mails: hushuangnian@163.com (S. Hu); sfhong@scu.edu.cn,\\
s-f.hong@tom.com, hongsf02@yahoo.com (S. Hong); qincn328@sina.com (X. Qin)}
\thanks{S. Hong is the corresponding author and was supported
partially by National Science Foundation of China Grant \#11371260.
X. Qin was supported partially by the Science and Technology Research Projects
of Chongqing Education Committee Grant \#KJ15012004.}}

\maketitle

{\bf Abstract.}
Let $\mathbb{F}_q$ denote the finite field of odd characteristic $p$ with $q$
elements ($q=p^{n},n\in \mathbb{N} $) and  $\mathbb{F}_q^*$
represent the nonzero elements of  $\mathbb{F}_{q}$. In this
paper, by using the Smith normal form we give an explicit formula
for the number of rational points of the algebraic variety defined
by the following system of equations over $\mathbb{F}_{q}$:
\begin{align*}
{\left\{\begin{array}{rl}
&\sum_{i=1}^{r_1}a_{1i}x_1^{e^{(1)}_{i1}} ...x_{n_1}^{e^{(1)}_{i,n_1}}
+\sum_{i=r_1+1}^{r_2}a_{1i}x_1^{e^{(1)}_{i1}} ...x_{n_2}^{e^{(1)}_{i,n_2}}-b_1=0,\\
&\sum_{j=1}^{r_3}a_{2j}x_1^{e^{(2)}_{j1}} ...x_{n_3}^{e^{(2)}_{j,n_3}}
+\sum_{j=r_3+1}^{r_4}a_{2j}x_1^{e^{(2)}_{j1}} ...x_{n_4}^{e^{(2)}_{j,n_4}}-b_2=0,
\end{array}\right.}
\end{align*}
where the integers $1\leq r_1<r_2$, $1\leq r_3<r_4$, $1\le n_1<n_2$, $1\le n_3<n_4$, $n_1\leq n_3$,
$b_1, b_2\in \mathbb{F}_{q}$, $a_{1i}\in \mathbb{F}_{q}^{*}$$(1\leq i\leq r_2)$,
$a_{2j}\in \mathbb{F}_{q}^{*}$$(1\leq j\leq r_4)$
and the exponent of each variable is positive integer. An example
is also presented to demonstrate the validity of the main result.

%%%%%%%%%%%%%%%%%%%%%%%%%%%%%%%%%%%%%%%%%%%%%%%%%%%%%%%%%%%



\section{Introduction}
Let $\mathbb{F}_{q}$ denote the finite field of $q$
elements with odd characteristic $p$ ($q=p^{n},n\in \mathbb{N} $
(the set of positive integers)) and $\mathbb{F}_{q}^*=\mathbb{F}_{q}\setminus\{0\}$.
Let $m$ be a positive integer, $f_i(x_1, ..., x_n)$$(i=1,...m)$
be some polynomial with $n$ variables over $\mathbb{F}_{q}$ and $V$ denote the
algebraic variety defined by the following system of equations over $\mathbb{F}_{q}$:
$$\left\{
\begin{aligned}
f_1(x_1, ..., x_n)&=0,\\
...&...\\
f_m(x_1, ..., x_n)&=0.
\end{aligned}
\right.
$$
By $N_q(V)$ we denote the number of $\mathbb{F}_q$-rational
points on the algebraic variety $V$ in $\mathbb{F}^n_q$. That is,
$$\mathbb{N}_q(V):=\#\{(x_1,...,x_n)
\in\mathbb{F}_{q}^n|f_i(x_1, ..., x_n)=0,\ i=1,...,m\}.$$
Particularly, if $m=1$, then we use $N_q(f)$ to stand for $N_q(V)$.

Studying the exact value of $N_q(V)$ is one of the main topics in finite fields.
The degrees $\text{deg}(f_i)$ play an important role in the estimate of $N_q(V)$.
Let $\lceil x\rceil$ denote the least integer
$\geq x$ and $\text{ord}_q$ denote the  additive valuation such that $\text{ord}_qq=1$.
In 1964, Ax \cite{[A]} generalized the Chevalley-Warning theroem by showing that
$$
\text{ord}_q N_q(V)\geq \Big\lceil\frac{n-\sum_{i=1}^{m}
\text{deg}f_i}{\sum_{i=1}^{m}\text{deg}f_i}\Big\rceil.
$$
Later on, further works were done by Katz \cite{[K]}, Adolphson-Sperber
\cite{[AS1]}-\cite{[AS2]}, Moreno-Moreno \cite{[MM]} and Wan \cite{[W1]}-\cite{[W3]}.

It is difficult to give an explicit formula for $N_q(V)$ in general.
Finding explicit formula for $N_q(f)$ under certain conditions
has attracted many authors for many years (see, for instance, \cite{[HV]} et al).
It is well known that there exists an explicit formula for $N_q(f)$ with
$\text{deg}(f) \leq2$ in $\mathbb{F}_{q}$ (see, for example, \cite{[IR]}
and \cite{[L]}). One first considered the diagonal hypersurface:
$$
a_1x_1^{e_1} + ... + a_nx_n^{e_n}-b=0, \ 1\leq i\leq n,
\ a_i\in \mathbb{F}_{q}^{*}, \ b\in \mathbb{F}_{q},\  e_i > 0,\eqno(1.1)
$$
and much work has been done to seek for the number of rational points of
the hypersurface (1.1)(see, for instance, \cite{[S1]} and \cite{[Wo1]}-\cite{[Wo2]}).
Carlitz \cite{[Ca1]}, Cohen \cite{[Co]} and Hodges \cite{[Hod]}
counted the rational points on the following $k$-linear hypersurface
$$
a_1x_{11}...x_{1k} + a_2x_{21}...x_{2k} + ... + a_nx_{n1}...x_{nk}-b=0,\eqno(1.2)
$$
with $a_i\in \mathbb{F}_{q}^{*}$$(1\leq i\leq n)$, $b\in \mathbb{F}_{q}$. Cao \cite{[C1]},
Cao and Sun \cite{[CaSu1]} \cite{[CaSu2]} studied the rational points
on the following more general diagonal hypersurface
$$
a_1x_{11}^{e_{11}}...x_{1n_1}^{e_{1n_1}} +a_2x_{21}^{e_{21}}...x_{2n_2}
^{e_{2n_2}}+ ... + a_rx_{r1}^{e_{r1}}...x_{rn_r}^{e_{rn_r}}=0,
$$
with $1\leq i\leq r$, $1\leq j\leq n_i$, $e_{ij}\in \mathbb{N}$,
$a_i\in \mathbb{F}_{q}^{*}$. Clearly, this extends (1.1) and (1.2)
when $b=0$. On the other hand, Baoulina \cite{[B1]} \cite{[B2]},
Pan, Zhao and Cao \cite{[PC]} considered the hypersurface of the form
$(a_1x_1^{m_1}+...+a_nx_n^{m_n})^\lambda -bx_1^{k_1}...x_n^{k_n}=0$
which extended the results of Carlitz in \cite{[Ca2]} and \cite{[Ca3]}.

If $f=a_1x_1^{e_{11}}...x_n^{e_{1n}} + ... + a_sx_1^{e_{s1}}...x_n^{e_{sn}}-b$
with $e_{ij}>0$ and $a_i\in \mathbb{F}_{q}^{*}$ for $1\leq i\leq s$ and $1\leq j\leq n$,
$b\in \mathbb{F}_{q}$, then a formula for $N_q(f)$ was given by Sun \cite{[S2]}.
Moreover, if $s=n$ and $\text{gcd}(\det(e_{ij}),q-1)=1$
($\det(e_{ij})$ represents the determinant of the $n\times n$ matrix
$\big(e_{ij}\big)$, then Sun \cite{[S2]} gave the explicit formula
for the number of rational points as follows:
\begin{align*}
N(f)={\left\{\begin{array}{rl}q^n-(q-1)^n+
\frac{(q-1)^n+(-1)^n(q-1)}{q},& \ \text{if} \ b=0,\\
\frac{(q-1)^n-(-1)^n}{q},  &  \ \text{if} \ b\neq0.
\end{array}\right.}
\end{align*}
Wang and Sun \cite{[WS]} gave the formula for the
number of rational points of the following hypersurface
$$
a_1x_1^{e_{11}} +a_2x_1^{e_{11}} x_2^{e_{22}}+ ... +
a_nx_1^{e_{n1}}x_2^{e_{n2}}...x_n^{e_{nn}}-b=0,
$$
with  $e_{ij}\geq0$, $a_i\in \mathbb{F}_{q}^{*}$, $b\in \mathbb{F}_{q}$.
In 2005, Wang and Sun \cite{[WS051]} extended  the results of
\cite{[S2]} and \cite{[WS]}. Recently, Hu, Hong and Zhao \cite{[HHZ]}
generalized Wang and Sun's results. In fact, they used the
Smith normal form to present a formula for $N_q(f)$ with $f$
being given by:
\begin{equation*}
f=\sum\limits_{j=0}^{t-1}\sum\limits_{i=1}^{r_{j+1}-r_j}
a_{r_j+i}x_1^{e_{r_j+i,1}} ...x_{n_{j+1}}^{e_{r_j+i,n_{j+1}}}-b, \eqno(1.3)
\end{equation*}
where the integers $t>0$, $r_0=0<r_1<r_2<...<r_t$, $1\le n_1<n_2<...<n_t$,
$b\in \mathbb{F}_{q}$, $a_i\in \mathbb{F}_{q}^{*}$ $(1\leq i\leq r_t)$
and the exponents $e_{ij}$ of each variable are positive integers.
On the other hand, Yang \cite{[Y]} followed Sun's method and gave a formula for
the rational points $N_q(V)$ on the following algebraic variety $V$ over $\mathbb{F}_{q}$:
$$\left\{
\begin{aligned}
a_{11}x_1^{e_{11}^{(1)}}...x_n^{e_{1n}^{(1)}} + ... +
a_{1s}x_1^{e_{s1}^{(1)}}...x_n^{e_{sn}^{(1)}}-b_1&=0,\\
...&...\\
a_{m1}x_1^{e_{11}^{(m)}}...x_n^{e_{1n}^{(m)}} + ...
+ a_{ms}x_1^{e_{s1}^{(m)}}...x_n^{e_{sn}^{(m)}}-b_m&=0.
\end{aligned}
\right.
$$
Song and Chen \cite{[SC]} continued to make use of
Sun's method and obtained a formula for $N_q(V)$ with $V$ being
the algebraic variety over $\mathbb{F}_{q}$ defined by:
$$\left\{
\begin{aligned}
\sum_{j=1}^{s_1}a_{1j}x_1^{e_{j1}^{(1)}}...x_{n_1}^{e_{j,{n_1}}^{(1)}}
+\sum_{j=s_1+1}^{s_2}a_{1j}x_1^{e_{j1}^{(1)}}...x_{n_2}^{e_{j,{n_2}}^{(1)}}-b_1&=0,\\
...&...\\
\sum_{j=1}^{s_1}a_{mj}x_1^{e_{j1}^{(m)}}...x_{n_1}^{e_{j,{n_1}}^{(m)}}
+\sum_{j=s_1+1}^{s_2}a_{mj}x_1^{e_{j1}^{(m)}}...x_{n_2}^{e_{j,{n_2}}^{(m)}}-b_m&=0.
\end{aligned}
\right.
$$
Meanwhile, they proposed an interesting question which was recently
answerd by Hu and Hong \cite{[HH]}. A more general question was
suggested by Hu, Hong and Zhao in \cite{[HHZ]} that can be stated
as follows.

{\bf Problem 1.1.} \cite{[HHZ]} Let $m\ge 1$ be an integer and $t_1, ..., t_m$
be positive integers. Find an explicit formula for the number of rational points
on the following algebraic variety over $\mathbb{F}_q$:
\begin{equation*}
\sum\limits_{j=0}^{t_k-1}\sum\limits_{i=1}^{r_{k,j+1}-r_{kj}}
a_{r_{kj}+i}x_1^{e^{(k)}_{r_{kj}+i,1}} ...x_{n_{k,j+1}}^{e^{(k)}
_{r_{kj}+i,n_{k,j+1}}}-b_k=0, \ k=1,...,m,
\end{equation*}
with $0=r_{k0}<r_{k1}<r_{k2}<...<r_{k,t_k}$, $1\le n_{k1}<n_{k2}<...<n_{k,t_k}$
and $e^{(k)}_{r_{kj}+i,j}>0$ being integers,
$b_k\in \mathbb{F}_{q}$, $a_{r_{kj}+i}\in \mathbb{F}_{q}^{*}$ for
$1\le k\le m$, $1\le i\le r_{k,t_k}$ and $1\leq j\leq t_k$.

If $m=1$, then the main result in \cite{[HHZ]} answers Problem 1.1. But it is
kept open when $m\ge 2$. Clearly, Yang \cite{[Y]}, Song and Chen \cite{[SC]} and Hu and
Hong \cite{[HH]} gave a partial answer to Problem 1.1 when $m\ge 2$.

In what follows, we always let $ r_i, n_i(i=1,...,4)$ be positive integers such that
$1\leq r_1<r_2$, $1\leq r_3<r_4$, $1\le n_1<n_2$, $1\le n_3<n_4$, $n_1\leq n_3$,
$b_1, b_2\in \mathbb{F}_{q}$, $a_{1i}\in \mathbb{F}_{q}^{*}$$(1\leq i\leq r_2)$
 and $a_{2j}\in \mathbb{F}_{q}^{*}$$(1\leq j\leq r_4)$. Let $n=\max\{n_2,n_4\}$,
$k=1,2$  and $f_k(\mathrm{x}):=f_k(x_1,...,x_{n})\in \mathbb{F}_q[x_1,...,x_{n}]$
be defined by
$$
f_1(\mathrm{x}):=f_1(x_1,...,x_{n})=\sum\limits_{i=1}^{r_2}a_{1i}
\mathrm{x}^{E^{(1)}_i}-b_1
$$
and
$$
f_2(\mathrm{x}):=f_2(x_1,...,x_{n})=\sum\limits_{j=1}^{r_4}a_{2j}
\mathrm{x}^{E^{(2)}_j}-b_2.            \eqno(1.3)
$$
with $E^{(1)}_i(1\leq i\leq r_2)$ and  $E^{(2)}_j(1\leq j\leq r_4)$
being the vectors of non-negative integer components of
dimension $n$ defined as follows:
$$\left\{
\begin{aligned}
E^{(1)}_1&=(e^{(1)}_{11},...,e^{(1)}_{1n_1},0,...,0),\  \
\mathrm{x}^{E^{(1)}_1}=x_1^{e^{(1)}_{11}}...x_{n_1}^{e^{(1)}_{1n_1}}, \\
......&...... \\
E^{(1)}_{r_1}&=(e^{(1)}_{r_1,1},...,e^{(1)}_{r_1,n_1},0,...,0),\  \
\mathrm{x}^{E^{(1)}_{r_1}}=x_1^{e^{(1)}_{r_1,1}}...x_{n_1}^{e^{(1)}_{r_1,n_1}}, \\
E^{(1)}_{r_1+1}&=(e^{(1)}_{r_1+1,1},...,e^{(1)}_{r_1+1,n_2},0,...,0),\  \
\mathrm{x}^{E^{(1)}_{r_1+1}}=x_1^{e^{(1)}_{r_1+1,1}}...x_{n_2}^{e^{(1)}_{r_1+1,n_2}},\\
......&...... \\
E^{(1)}_{r_2}&=(e^{(1)}_{r_2,1},...,e^{(1)}_{r_2,n_2},0,...,0), \  \
\mathrm{x}^{E^{(1)}_{r_2}}=x_1^{e^{(1)}_{r_2,1}}...x_{n_2}^{e^{(1)}_{r_2,n_2}},\\
E^{(2)}_1&=(e^{(2)}_{11},...,e^{(2)}_{1n_3},0,...,0),\  \
\mathrm{x}^{E^{(2)}_1}=x_1^{e^{(2)}_{11}}...x_{n_3}^{e^{(2)}_{1n_3}}, \\
......&...... \\
E^{(2)}_{r_3}&=(e^{(2)}_{r_3,1},...,e^{(2)}_{r_3,n_3},0,...,0),\  \
\mathrm{x}^{E^{(2)}_{r_3}}=x_1^{e^{(2)}_{r_3,1}}...x_{n_3}^{e^{(2)}_{r_3,n_3}}, \\
E^{(2)}_{r_3+1}&=(e^{(2)}_{r_3+1,1},...,e^{(2)}_{r_3+1,n_4},0,...,0),\  \
\mathrm{x}^{E^{(2)}_{r_3+1}}=x_1^{e^{(2)}_{r_3+1,1}}...x_{n_4}^{e^{(2)}_{r_3+1,n_4}},\\
......&...... \\
E^{(2)}_{r_4}&=(e^{(2)}_{r_4,1},...,e^{(2)}_{r_4,n_4},0,...,0), \  \
\mathrm{x}^{E^{(2)}_{r_4}}=x_1^{e^{(2)}_{r_4,1}}...x_{n_4}^{e^{(2)}_{r_4,n_4}}.\\
\end{aligned}
\right.
$$
Then a special case of Problem 1.1 is the following question.

{\bf Problem 1.2.} Let $f_1(\mathrm{x})$ and $f_2(\mathrm{x})$
be given as in (1.3). What is the formula for the number of rational points
on the algebraic variety defined by the following system of equations over $\mathbb{F}_{q}$:

$$\left\{
\begin{aligned}
f_1(\mathrm{x})&=0,\\
f_2(\mathrm{\mathrm{x}})&=0?
\end{aligned}
\right.\eqno(1.4)
$$

When $n_1=n_3$ and $n_2=n_4$, this question was answered
by Song and Chen \cite{[SC]}. However, if $n_1\neq n_3$
or $n_2\neq n_4$, then Problem 1.2 has not been solved yet.

In this paper, our main goal is to investigate Problem 1.2.
We will follow and develop the method of \cite{[HHZ]} to
study Problem 1.2. To state the main result, we first
need to introduce some related concept and notation.
Throughout this paper, we let
\begin{equation*}
E_{11}:=
\left( \begin{array}{*{15}c}
E^{(1)}_1\\
\vdots\\
E^{(1)}_{r_1}\\
\end{array}\right)_{r_1\times m},
\ \ \ E_{12}:=
\left( \begin{array}{*{15}c}
E^{(1)}_{r_1+1}\\
\vdots\\
E^{(1)}_{r_2}\\
\end{array}\right)_{(r_2-r_1)\times m},
\end{equation*}

\begin{equation*}
E_{21}:=
\left( \begin{array}{*{15}c}
E^{(2)}_1\\
\vdots\\
E^{(2)}_{r_3}\\
\end{array}\right)_{r_3\times m}
\ {\rm and} \ E_{22}:=
\left( \begin{array}{*{15}c}
E^{(2)}_{r_3+1}\\
\vdots\\
E^{(2)}_{r_4}\\
\end{array}\right)_{(r_4-r_3)\times m}.
\end{equation*}
For any integer $l$ with $1\leq l\leq 4$, we let $E^{(n_l)}$
denote the remaining of the exponent matrix of (1.4)
which deletes the items containing the variable $x_{n_l+1}$.
In fact, one has
$$
E^{(n_1)}:=E_{11}\ {\rm if}\ n_1<n_3,
\ {\rm and} \ \left( \begin{array}{*{15}c}
E_{11}\\
E_{21}
\end{array}\right)
\ {\rm if}\ n_1=n_3,
$$
$$
E^{(n_3)}:=\left( \begin{array}{*{15}c}
E_{11}\\
E_{21}
\end{array}\right)\ {\rm if}\ n_3<n_2,
\ {\rm and} \ \left( \begin{array}{*{15}c}
E_{11}\\
E_{12}\\
E_{21}
\end{array}\right)
\ {\rm otherwise},
$$
$$
E^{(n_4)}:=\left( \begin{array}{*{15}c}
E_{11}\\
E_{21}\\
E_{22}\\
\end{array}\right)\ {\rm if}\ n_4<n_2,
\ {\rm and} \ \left( \begin{array}{*{15}c}
E_{11}\\
E_{12}\\
E_{21}\\
E_{22}
\end{array}\right)
\ {\rm otherwise}.
$$
If $n_3\leq n_2$, then
\begin{equation*}
E^{(n_2)}:=
\left( \begin{array}{*{15}c}
E_{11}\\
E_{12}\\
E_{21}
\end{array}\right)
\ {\rm if} \ n_2<n_4,
\ {\rm and} \
\left( \begin{array}{*{15}c}
E_{11}\\
E_{12}\\
E_{21}\\
E_{22}
\end{array}\right)
\ {\rm otherwise},
\end{equation*}
and if $n_3>n_2$, then
$$E^{(n_2)}:=\left( \begin{array}{*{15}c}
E_{11}\\
E_{12}
\end{array}\right).$$
Then the Smith normal form (see \cite{[S]}, \cite{[Hu]}
or Section 2 below) guarantees the existences of unimodular matrices
$U^{(n_l)}$ and $V^{(n_l)}$  such that
\begin{equation*}
U^{(n_l)}E^{(n_l)}V^{(n_l)}=
\left( \begin{array}{*{15}c}
D^{(n_l)}&0\\
0&0
\end{array}\right),
\end{equation*}
where $D^{(n_l)}:={\rm diag}(d^{(l)}_{1}, ..., d^{(l)}_{s_l})$
with all the diagonal elements $d^{(l)}_{1}, ..., d^{(l)}_{s_l}$
being positive integers such that $d^{(l)}_{1}|...|d^{(l)}_{s_l}$.
Throughout, pick $\alpha \in\mathbb{F}^*_q$ to be a fixed primitive
element of $\mathbb{F}_q$. For any $\beta\in\mathbb{F}^*_q$, there
exists exactly an integer $r\in [1, q-1]$ such that $\beta=\alpha^r$.
Such an integer $r$ is called {\it index} of $\beta$ with
respect to the primitive element $\alpha $, and is denoted by
$\text{ind}_\alpha(\beta):=r$.

For $l=1,...,4$, we let $SE^{(n_l)}$ denote the system (1.4)
which deletes the items containing the variable $x_{n_l+1}$.
Clearly the exponent matrix of $SE^{(n_l)}$ is $E^{(n_l)}$.
For example, if $l=1$ and $n_1<n_3$, then $SE^{(n_1)}$
becomes the following system of equations over $\mathbb{F}_{q}$:
\begin{align*}
{\left\{\begin{array}{rl}
&\displaystyle\sum_{i=1}^{r_1}a_{1i}x_1^{e^{(1)}_{i1}} ...x_{n_1}^{e^{(1)}_{i,n_1}}-b_1=0,\\
&b_2=0,
\end{array}\right.}
\end{align*}
where $a_{1i}$, $b_1$, $b_2$ and $e^{(1)}_{ij}$ are given 
in (1.4) for all integers $i$ and $j$ with
$1\leq i\leq r_1$ and $1\leq j\leq n_1$.
Let $x_1^{e^{(1)}_{i1}} ...x_{n_1}^{e^{(1)}_{i,n_1}}=u_{1i}$ 
for any integer $i$ with $1\leq i\leq r_1$.
If $b_2\neq 0$, then it is obvious that $SE^{(n_1)}$ has 
no solution. If $b_2=0$, then $SE^{(n_1)}$ becomes
$$
\displaystyle\sum_{i=1}^{r_1}a_{1i}u_{1i}-b_1=0.\eqno (\overline{SE}^{(n_1)})
$$
Lemma 2.4 tells us the formula for the number of the solution
$(u_{11},...,u_{1r_1})\in ({\mathbb{F}}_{q}^{*})^{r_1}$ of $(\overline{SE}^{(n_1)})$.
Now we use $H_1$  to denote the number of the solutions
$(u_{11},...,u_{1r_1})\in ({\mathbb{F}}_{q}^{*})^{r_1}$ of $(\overline{SE}^{(n_1)})$
under the following extra conditions:
$$
{\left\{\begin{array}{rl}\text{gcd}(q-1,d_i^{(1)})|h_i'& {\rm for} \ i=1,...,s_1,\\
(q-1)|h_i'&{\rm for} \  i=s_1+1,...,r_1,
\end{array}\right.}\eqno (\widetilde{SE}^{(n_1)})
$$
where $(h'_1, ...,h'_{r_1})^T=U^{(n_1)}(\text{ind}_\alpha(u_{11}),...,
\text{ind}_\alpha(u_{1r_1}))^T$. Similarly, for $l=2,3,4$,
we can get the system of equations $SE^{(n_l)}$ and use $H_l$ to
denote the number of the solutions of $(\overline{SE}^{(n_l)})$
under the extra conditions $(\widetilde{SE}^{(n_l)})$.
One knows that $H_l$ (see \cite{[HH]}) is independent
of the choice of the primitive element $\alpha$.
In what follows, we let $N$ stand for the number of
rational points on the algebraic variety defined by (1.4).

Now let
$$N_0:=q^{\max\{n_2,n_4\}-n_1}(q^{n_1}-(q-1)^{n_1}) \ {\rm if} \
b_1=b_2=0, \ {\rm and} \ 0 \ {\rm otherwise},$$
\begin{align*}
N_{1}:={\left\{\begin{array}{rl}
q^{\max\{n_2,n_4\}-\min\{n_2,n_3\}}(q^{\min\{n_2,n_3\}-n_1}
& -(q-1)^{\min\{n_2, n_3\}-n_1})L_1,\\
&  {\rm if} \ n_1<n_3 \ {\rm and}\ b_2=0,\\
q^{\max\{n_2,n_4\}-\min\{n_2,n_4\}}(q^{\min\{n_2,n_4\}-n_1}
&-(q-1)^{\min\{n_2,n_4\}-n_1})L_1,\\
& {\rm if} \ n_1=n_3 \ {\rm and}\ b_2=0, \\
0,  & {\rm otherwise},
\end{array}\right.}
\end{align*}
$$
N_2:=q^{\max\{n_2,n_4\}-\min\{n_2,n_4\}}
(q^{\min\{n_2,n_4\}-n_3}-(q-1)^{\min\{n_2,n_4\}-n_3})L_3$$
$$
\ {\rm if} \ n_3<n_2, \ {\rm and} \ 0 \ {\rm otherwise},
$$
\begin{align*}
N_3:={\left\{\begin{array}{rl}\big(q^{n_4-n_2}-(q-1)^{n_4-n_2}\big)L_2, 
&  {\rm if} \ n_1\leq n_3<n_2<n_4,\\
\big(q^{n_4-n_3}-(q-1)^{n_4-n_3}\big)L_3, &  {\rm if} \ n_3\geq n_2, \\
0,  & {\rm otherwise},
\end{array}\right.}
\end{align*}
$N_4:=L_4$ if $n_4\geq n_2$, and $L_2$ otherwise, 
$N_5:=\big(q^{n_2-n_4}-(q-1)^{n_2-n_4}\big)L_4$
if $n_2>n_4$, and $0$ otherwise, $N_6:=q^{n_4-n_3}
(q^{n_3-n_2}-(q-1)^{n_3-n_2})L_2$ if $n_3>n_2$,
and $b_2=0$, and $0$ otherwise, where  for $1\le i\le 4$,
$$L_i:=H_i(q-1)^{n_i-s_i}\prod\limits_{j=1}^{s_i}
\gcd(q-1,d_j^{(i)}).$$

We can now state the main result of this paper 
which answers Problem 1.2 completely.

{\bf Theorem 1.3.} {\it Let $N$ denote the number of rational
points on the algebraic variety (1.4). Then
\begin{align*}
N={\left\{\begin{array}{rl}
\sum\limits_{i=0}^{4}N_i,& \  if  \ n_1<n_3<n_2<n_4,\\
\sum\limits_{i=0\atop i\neq2}^4N_i, & \  if  \ n_2=n_3,\\
\sum\limits_{i=0\atop i\neq1}^4N_i,& \  if  \ n_1=n_3,\ n_2<n_4,\\
\sum\limits_{i=0\atop i\neq3}^4N_i,& \  if  \ n_1<n_3,\ n_2=n_4,\\
\sum\limits_{i=0\atop i\neq1, 3}^4N_i, & \  if  \ n_1=n_3,\ n_2=n_4,\\
\sum\limits_{i=0\atop i\neq3}^{5}N_i,& \  if  \ n_1<n_3,\ n_2>n_4,\\
\sum\limits_{i=0\atop i\neq1,3}^{5}N_i,& \  if  \ n_1=n_3,\ n_2>n_4,\\
\sum\limits_{i=0\atop i\neq2,5}^{6}N_i,& \  if  \ n_3>n_2.
\end{array}\right.}
\end{align*}
with $N_{i}(0\le i\le 6)$ being defined as above.}

We organize this paper as follows. In Section 2, we present some useful
lemmas which will be needed later. In fact, we will first recall some
basic facts on the Smith normal form of an integer matrix.
Then we can use them to give a formula for the number of the system of
linear congruences with the same modulo. Consequently, in Section 3,
we first show a key lemma and then use it to prove Theorem 1.3.
In the final section, we provide an example to demonstrate
the validity of Theorem 1.3.

Throughout this paper, $\gcd(a,m)$ will denote the greatest
common divisor of any positive integers $a$ and $m$.

%%%%%%%%%%%%%%%%%%%%%%%%%%%%%%%%%%%%%%%%%%%%%%%%%%%%%%%%%%%%%%%%%%%%%%%%%%%%%%%%%%
\section{Preliminary lemmas}
In this section, we present some useful lemmas that are needed in Section 3.
We first recall two well-known definitions.

{\bf Definition 2.1.} \cite{[Hu]} Let $M$ be a square integer matrix. If
the determinant of $M$ is $\pm 1$, then $M$ is called an {\it unimodular matrix}.

{\bf Definition 2.2.} \cite{[Hu]}  For given any positive integers $m$ and $n$, let
$P$ and $Q$ be two $m\times n$ integer matrices. Suppose that there are two
modular matrices $U$ of order $m$ and  $V$ of order $n$ such that
$P=UQV$. Then we say that $P$ and $Q$ are {\it equivalent}
and we write $P\thicksim Q$.

Clearly equivalence has the three properties of being
reflexive, symmetric and transitive.

{\bf Lemma 2.1.} \cite{[Hu]} \cite{[S]}
{\it Let $P$ be a nonzero $m\times n$ integer matrix.
Then $P$ is equivalent to a block matrix of the following form
\begin{equation*}
\left( \begin{array}{*{15}c}
D&0\\
0&0
\end{array}\right), \eqno(2.1)
\end{equation*}
where $D={\rm diag}(d_1, ..., d_r)$ with all the diagonal elements $d_i$
being positive integers and satisfying that $d_i|d_{i+1}$ $(1 \le i < r)$.
In other words, there are unimodular matrices $U$ of order $m$ and $V$
of order $n$ such that}
$$
UPV=\left( \begin{array}{*{15}c}
{\rm diag}(d_1, ..., d_r)&0\\
0&0
\end{array}\right).
$$

We call the diagonal matrix in (2.1) the {\it Smith normal form} of the matrix $P$.
Usually, one writes the Smith normal form of $P$ as ${\rm SNF}(P)$.
The elements $d_i$ are unique up to multiplication by a unit and are called
the {\it elementary divisors, invariants,} or {\it invariant factors}.

For any system of linear congruences
$$
{\left\{\begin{array}{rl}
\sum_{j=1}^{n}h_{1j}y_j &\equiv b_1\pmod m,\\
............\\
\sum_{j=1}^{n}h_{sj}y_j &\equiv b_s\pmod m,
\end{array}\right.}  \eqno (2.2)
$$
let $Y=(y_1,..., y_n)^T$ be the column of indeterminates $y_1, ..., y_n$,
$B=(b_1,...,b_s)^T$ and $H=(h_{ij})$ be the matrix of its
coefficient. Then one can write (2.2) as
$$HY\equiv B\pmod m. \eqno(2.3)$$
By Lemma 2.1, there are unimodular matrices $U$ of order $s$ and $V$
of order $n$ satisfying that
$$
UHV={\rm SNF}(H)=\left( \begin{array}{*{15}c}
{\rm diag}(d_1, ..., d_r)&0\\
0&0
\end{array}\right).
$$

Finally, we present the following two known lemmas.

{\bf Lemma 2.2.} \cite{[HHZ]}
{\it Let $B'=(b'_1,...,b'_s)^T=UB$. Then the system (2.3)
of linear congruences is solvable if and only if
$\gcd(m, d_i)|b'_i$ for all integers $i$ with $1\le i\le r$ and $m|b'_i$
for all integers $i$ with $r+1\le i\le s$. Further,
the number of solutions of (2.3) is equal to
$m^{n-r}\prod_{i=1}^r\gcd(m, d_i).$}

{\bf Lemma 2.3.} \cite{[S2]}
{\it Let $c_1,...,c_k\in \mathbb{F}_q^{*} $ and  $c \in\mathbb{F}_{q}$, and
let $N(c)$ denote the number of rational points $(x_1, ..., x_k)\in
(\mathbb{F}_{q}^{*})^k$ on the hypersurface $c_1x_1+...+c_kx_k=c.$ Then}
\begin{align*}
N(c)={\left\{\begin{array}{rl}\frac{(q-1)^k+(-1)^k(q-1)}{q},& \ {\it if} \ c=0,\\
\frac{(q-1)^k-(-1)^k}{q},  & {\it otherwise}.
\end{array}\right.}
\end{align*}

{\bf Lemma 2.4.}
{\it Let $c_{ij}\in \mathbb{F}_q^{*}$ for all integers $i$ and $j$ with
$1\le i\le m$ and $1\le j\le k$ and  $c_1, ..., c_m\in\mathbb{F}_{q}$.
Let $N(c_1, ..., c_m)$ denote the number of rational points
$(x_{11}, ..., x_{1k}, ..., x_{m1}, ..., x_{mk})\in (\mathbb{F}_{q}^{*})^{mk}$
on the following variety
$$
{\left\{\begin{array}{rl}
&c_{11}x_{11}+...+c_{1k}x_{1k}=c_1,\\
&......\\
&c_{m1}x_{m1}+...+c_{mk}x_{mk}=c_m.
\end{array}\right.} \eqno(2.4)
$$
Then
$$
N(c_1, ..., c_m)=\frac{(q-1)^r}{q^m}((q-1)^{k-1}+(-1)^k)^r((q-1)^k-(-1)^k)^{m-r},
$$
where $r:=\#\{1\le i\le m| c_i=0\}$.}

{\it Proof.} For $1\le i\le m$, let $N(c_i)$ denote the number of
rational points $(x_{i1}, ..., x_{ik})\in (\mathbb{F}_{q}^{*})^{k}$
on the hypersurface $c_{i1}x_{i1}+...+c_{ik}x_{ik}=c_i$. Since
for any rational points
$(x_{11}, ..., x_{1k}, ..., x_{m1}, ..., x_{mk})\in (\mathbb{F}_{q}^{*})^{mk}$
on the variety (2.4), the different part of coordinates
$(x_{i1}, ..., x_{ik})\in (\mathbb{F}_{q}^{*})^{k} (1\le i\le m)$
are independent, one has
$$N(c_1, ..., c_m)=\prod_{i=1}^mN(c_i).$$
So Lemma 2.3 applied to $N(c_i)$ gives us the required result.
Hence Lemma 2.4 is proved. \hfill$\Box$

%%%%%%%%%%%%%%%%%%%%%%%%%%%%%%%%%%%%%%%%%%%%%%%%%%%%%%%%%%%%%%%%%%%%%%%%%%%%%%%%%%
\section{Proof of Theorem 1.3}

In this section, we show Theorem 1.3. First, we provide
some notation and two lemmas. For any given $(u_{11},...,u_{1r_2},...,u_{21},
...,u_{2r_4})\in \mathbb{ F}_{q}^{r_2+r_4}$, we use
$$
N(\mathrm{x}^{E^{(1)}_i}=u_{1i}(i=1,2,..., r_2),\mathrm{x}^{E^{(2)}_j}=u_{2j}(j=1,2,..., r_4 ).)
$$
to denote the number of rational points
$(x_1,...,x_{\max\{n_2,n_4\}})\in\mathbb{ F}_{q}^{\max\{n_2,n_4\}}$
of the following algebraic variety over $\mathbb{ F}_{q}$:
$$\left\{
\begin{aligned}
\mathrm{x}^{E^{(1)}_1}&=u_{11}, \\
.......&..... \\
\mathrm{x}^{E^{(1)}_{r_2}}&=u_{1r_2}, \\
\mathrm{x}^{E^{(2)}_1}&=u_{21},\\
......&...... \\
\mathrm{x}^{E^{(2)}_{r_4}}&=u_{2r_4}.
\end{aligned}
\right. \eqno(3.1)
$$
It follows that
$$
N=\sum_{\begin{subarray}{I}(u_{11},...,u_{1r_2},u_{21},...,u_{2r_4})\in \mathbb{ F}_{q}^{r_2+r_4},\\
a_{11}u_{11}+...+a_{1r_2}u_{1r_2}=b_1,\\
a_{21}u_{21}+...+a_{2r_4}u_{2r_4}=b_2.
\end{subarray}}N\Big(\left\{
\begin{aligned} \mathrm{x}^{E^{(1)}_i}=u_{1i}(i=1,2,..., r_2),\\
\mathrm{x}^{E^{(2)}_j}=u_{2j}(j=1,2,...,r_4).
\end{aligned}
\right.\Big).
$$
Define
$$T:=\{(u_{11},...,u_{1r_2},u_{21},...,u_{2r_4})\in \mathbb{ F}_{q}^{r_2+r_4}:\begin{subarray}{I}
a_{11}u_{11}+...+a_{1r_2}u_{1r_2}=b_1,\\
a_{21}u_{21}+...+a_{2r_4}u_{2r_4}=b_2.\end{subarray}\},$$
with $b_1$, $b_2$, $a_{1i}$ $(1\leqslant i\leqslant r_2)$ and $a_{2j} (1\leqslant j\leqslant r_4)$
being given as in (1.4). Then
$$
N=\sum_{(u_{11},...,u_{1r_2},u_{21},...,u_{2r_4})\in T}N\Big(\left\{
\begin{aligned} \mathrm{x}^{E^{(1)}_i}=u_{1i}(i=1,2,..., r_2),\\
\mathrm{x}^{E^{(2)}_j}=u_{2j}(j=1,2,...,r_4).
\end{aligned}
\right.\Big).
\eqno(3.2)
$$
Let $(u_{11},...,u_{1r_2},u_{21},...,u_{2r_4})\in T$. We consider the algebraic variety of (3.1).
For integer $i$, $j$ with $1\leq i<r_2$ and $1\leq j<r_4$, one can easily deduce that
\begin{equation*}
u_{1,i+1}=0 \ {\rm if} \ u_{1i}=0, \ {\rm and} \ u_{2,j+1}=0 \ {\rm if}\ u_{2j}=0. \eqno(3.3)
\end{equation*}
Define
$$
V_1:=\{(u_{11},...,u_{1r_1},0,...,0)
\in\mathbb{F}_q^{r_2+r_4}: u_{1i}\in\mathbb{F}_q^*,1\leq i\leq r_1\},
$$
$$
V_2:=\{(u_{11},...,u_{1r_1},\underbrace{0,...,0}_{r_2-r_1},u_{21},...,u_{2r_3},0,...,0)
\in\mathbb{F}_q^{r_2+r_4}: u_{1i},u_{2j}\in\mathbb{F}_q^*,1\leq i\leq r_1,1\leq j\leq r_3\},
$$
$$
V_3:=\{(u_{11},...,u_{1r_2},u_{21}, ..., u_{2r_3}, 0,...,0)
\in\mathbb{F}_q^{r_2+r_4}: u_{1i},u_{2j}\in\mathbb{F}_q^*, 1\leq i\leq r_2, 1\leq j\leq r_3\},
$$
$$
V_4:=\{(u_{11},...,u_{1r_2},u_{21},...,u_{2r_4})
\in\mathbb{F}_q^{r_2+r_4}: u_{1i}, u_{2j}\in\mathbb{F}_q^*, 1\leq i\leq r_2, 1\leq j\leq r_4\},
$$
$$
V_5:=\{(u_{11},...,u_{1r_1},0,...,0,u_{21},...,u_{2r_4})
\in\mathbb{F}_q^{r_2+r_4}: u_{1i},u_{2j}\in\mathbb{F}_q^*,1\leq i\leq r_1,1\leq j\leq r_4\},
$$
and
$$
V_6:=\{(u_{11},...,u_{1r_2},0,...,0)
\in\mathbb{F}_q^{r_2+r_4}: u_{1i}\in\mathbb{F}_q^*,1\leq i\leq r_2\}.
$$
Let $T(V_0)$ consist of zero vector of dimension $r_2+r_4$
and $T(V_i)$ ($0\le i\le 6$) denote the set of the vectors $v$ such that
$v\in T$ and $v\in V_i$. For any integer $i$ with $0\le i\le 6$, define
$$M_i:=\sum\limits_{(u_{11}, ..., u_{1r_2}, u_{21}, ..., u_{2r_4})
\in T(V_i)}N\Big(\left\{
\begin{aligned} \mathrm{x}^{E^{(1)}_i}=u_{1i}(i=1,2,..., r_2),\\
\mathrm{x}^{E^{(2)}_j}=u_{2j}(j=1,2,...,r_4).
\end{aligned}
\right.\Big).\eqno(3.4)
$$
Then we have the following key lemma.

{\bf Lemma 3.1.} {\it One has that $M_i=N_i\ \  for\  all\  0\le i\le 6.$}

{\it Proof.} First we show that $M_0=N_0$.
For any given integer $i$ with $i=1,2$,
if there exists $b_i\neq0$, then the algebraic variety
\begin{align*}
{\left\{\begin{array}{rl}
&a_{11}u_{11}+...+a_{1r_2}u_{1r_2}-b_1=0,\\
&a_{21}u_{21}+...+a_{2r_4}u_{2r_4}-b_2=0,
\end{array}\right.}
\end{align*}
has no zero solutions $(u_{11},...,u_{1r_2},u_{21},...,u_{2r_4})
\in \mathbb{ F}_{q}^{r_2+r_4}$. So $T(V_0)$ is empty.
It follows that $M_{0}=0=N_0$. If $b_1=b_2=0$, then $T(V_0)$
consists of zero vector of dimension $r_2+r_4$. Then
\begin{align*}
M_0&=N\Big(\left\{
\begin{aligned} \mathrm{x}^{E^{(1)}_i}=u_{1i}(i=1,..., r_2)\\
\mathrm{x}^{E^{(2)}_j}=u_{2j}(j=1,...,r_4)
\end{aligned}
\right.\Big)\\
&=N(\mathrm{x}^{E^{(1)}_1}=0) \ \ (\text{since} \ n_1\leq n_3)\\
&=N(x_1^{e^{(1)}_{11}}...x_{n_1}^{e^{(1)}_{1n_1}}=0)\\
&=q^{\max\{n_2,n_4\}-n_1}\sum_{j=1}^{n_1}\binom{n_1}{j}(q-1)^{n_1-j}\\
&=q^{\max\{n_2,n_4\}-n_1}(q^{n_1}-(q-1)^{n_1})=N_0.
\end{align*}
This completes the proof of the first part of the Lemma 3.1.

We can now turn our attention to prove that $M_1=N_1$. First we we consider the case
$n_1<n_3<n_2$. Since $n_1<n_3$, it follows from (3.4) and the definition of $T(V_1)$ that
\begin{align*}
M_1=\sum_{\begin{subarray}{I}(u_{11},...,u_{1r_1})\in (\mathbb{F}_q^*)^{r_1},\\
a_{11}u_{11}+...+a_{1r_1}u_{1r_1}=b_1\ {\rm and}\ 0=b_2.
\end{subarray}}
N\Big(\left\{
\begin{aligned} \mathrm{x}^{E^{(1)}_i}&=u_{1i},\ 1\leq i\leq r_1,\\
\mathrm{x}^{E^{(1)}_j}&=0,\ r_1+1\leq j\leq r_2,\\
\mathrm{x}^{E^{(2)}_k}&=0,\ 1\leq k\leq r_4.
\end{aligned}
\right.\Big).
 \ \ (3.5)
\end{align*}

From (3.5), we deduce that if $b_2\neq0$, then $M_1=0=N_1$. That is,
$$M_1=N_1 \ {\rm if} \ n_1<n_3<n_2 \ {\rm and} \ b_2\neq0. $$

Now we let $b_2=0$. Notice that $n_3<n_2$,  
then the definitions of $\mathrm{x}^{E^{(1)}_j}$
and $\mathrm{x}^{E^{(2)}_k}$ tell us that
the fact that $\mathrm{x}^{E^{(1)}_j}=0$ and 
$\mathrm{x}^{E^{(2)}_k}=0$ with $r_1+1\le j\le r_2$, 
$1\leq k\leq r_4$ is reduced to saying that 
$\mathrm{x}^{E^{(2)}_1}=0$. It then follows from (3.5) that
\begin{align*}
&M_1=\sum_{\begin{subarray}{I}(u_{11},...,u_{1r_1})\in (\mathbb{F}_q^*)^{r_1},\\
a_{11}u_{11}+...+a_{1r_1}u_{1r_1}=b_1.
\end{subarray}}
N\Big(\mathrm{x}^{E^{(1)}_i}=u_{1i},\ 1\leq i\leq r_1\ {\rm and}\ \mathrm{x}^{E^{(2)}_1}=0\Big).
 \ \ (3.6)
\end{align*}

It is easy to see that $\mathrm{x}^{E^{(2)}_1}=0$ is equivalent to $x_1...x_{n_3}=0$.
Since $u_{1r_1}\ne 0$ and $u_{1r_1}= x_1^{e^{(1)}_{r_1,1}}...x_{n_1}
^{e^{(1)}_{r_1,n_1}}$, one has $x_1...x_{n_1}\ne 0$.
So $\mathrm{x}^{E^{(2)}_1}=0$ is equivalent to $x_{n_1+1}...x_{n_3}=0$.
Then by (3.6), one gets that
\begin{align*}
M_1=\sum_{\begin{subarray}{I}(u_{11},...,u_{1r_1})\in (\mathbb{F}_q^*)^{r_1},\\
a_{11}u_{11}+...+a_{1r_1}u_{1r_1}=b_1.
\end{subarray}}
N\Big(\mathrm{x}^{E^{(1)}_i}=u_{1i},\ 1\leq i\leq r_1\ {\rm and}\ x_{n_1+1}...x_{n_3}=0\Big).
 \ \ (3.7)
\end{align*}

For any given $(u_{11},...,u_{1r_1})\in (\mathbb{ F}_{q}^*)^{r_1}$ 
with $\sum_{i=1}^{r_1}a_{1i}u_{1i}=b_1$, one has
\begin{align*}
&N(\mathrm{x}^{E^{(1)}_i}=u_{1i}, 1\leq i\leq r_1\ {\rm and} \ x_{n_1+1}...x_{n_3}=0)\\
&=\#\{(x_1,...,x_{\max\{n_2,n_4\}})\in (\mathbb{F}_q)^{\max\{n_2,n_4\}}: 
\mathrm{x}^{E^{(1)}_i}=u_{1i}, 1\leq i\leq r_1\ {\rm and} \ x_{n_1+1}...x_{n_3}=0\}.
\end{align*}
Since each of the components $x_{n_3+1},..., x_{\max\{n_2,n_4\}}$ can run over
the whole finite field $\mathbb{F}_q$ independently, it then follows that
\begin{align*}
&N(\mathrm{x}^{E^{(1)}_i}=u_{1i}, 1\leq i\leq r_1\ {\rm and} \ x_{n_1+1}...x_{n_3}=0)\\
&=q^{\max\{n_2,n_4\}-n_3}\times \\
&\#\{(x_1,...,x_{n_3})\in (\mathbb{F}_q)^{n_3}: \mathrm{x}^{E^{(1)}_i}=u_{1i}, 1\leq i\leq r_1\
{\rm and} \ x_{n_1+1}...x_{n_3}=0\}.\ \ \ \ \ (3.8)
\end{align*}
Notice that the choice of $(x_1,...,x_{n_1})\in (\mathbb{F}_q^*)^{n_{1}}$
satisfying that $\mathrm{x}^{E^{(1)}_i}=u_{1i}$ ($i=1,...,r_1$)
is independent of the choice of
$(x_{n_1+1},...,x_{n_3})\in (\mathbb{F}_q)^{n_3-n_1}$
satisfying that $x_{n_1+1}...x_{n_3}=0$. We then derive that
\begin{align*}
&\#\{(x_1,...,x_{n_3})\in (\mathbb{F}_q)^{n_3}: \mathrm{x}^{E^{(1)}_i}=u_{1i},1\leq i\leq r_1\ {\rm and} \ x_{n_1+1}...x_{n_3}=0\}\\
=&\#\{(x_1,...,x_{n_1})\in (\mathbb{F}_q^*)^{n_1}: \mathrm{x}^{E^{(1)}_i}=u_{1i}, 1\leq i\leq r_1\}\times \\
&\#\{(x_{n_1+1},...,x_{n_3})\in (\mathbb{F}_q)^{n_3-n_1}: x_{n_1+1}...x_{n_3}=0\}.
\ \ \ \ \ \ \ \ \ \ \ \ \ \ \ \ \ \ \ \ \ \ \ \ \ \ \  \ \ \ (3.9)
\end{align*}

On the other hand, we can easily compute that
\begin{align*}
&\#\{(x_{n_1+1},...,x_{n_3})\in (\mathbb{F}_q)^{n_3-n_1}: x_{n_1+1}...x_{n_3}=0\}\\
&=\sum_{i=1}^{n_3-n_1}{n_3-n_1 \choose i}(q-1)^{n_3-n_1-i}=q^{n_3-n_1}-(q-1)^{n_3-n_1}.
\ \ \ \ \ \ \ \ \ \ \ \ \ \ \ \ \ \ (3.10)
\end{align*}
So by (3.8) to (3.10), one obtains that
\begin{align*}
&N(\mathrm{x}^{E^{(1)}_i}=u_{1i}, 1\leq i\leq r_1\ {\rm and} \ x_{n_1+1}...x_{n_3}=0)\\
&=q^{\max\{n_2,n_4\}-n_3}(q^{n_3-n_1}-(q-1)^{n_3-n_1})\\
&\times \#\{(x_1,...,x_{n_{1}})\in (\mathbb{F}_q^*)^{n_1}: \mathrm{x}^{E^{(1)}_i}=u_{1i}, 1\leq i\leq r_1\}\\
&=q^{\max\{n_2,n_4\}-n_3}(q^{n_3-n_1}-(q-1)^{n_3-n_1}) N(\mathrm{x}^{E^{(1)}_i}=u_{1i}, 1\leq i\leq r_1). \ \ \ \ \ \ \ \ \ \ \  (3.11)
\end{align*}
Then by (3.7) together with (3.11), we have
\begin{align*}
M_1=&q^{\max\{n_2,n_4\}-n_3}(q^{n_3-n_1}-(q-1)^{n_3-n_1})\times \\
&\sum_{\begin{subarray}{I}(u_{11},...,u_{1r_1})\in (\mathbb{ F}_{q}^*)^{r_1},\\
a_{11}u_{11}+...+a_{{1r_1}}u_{1r_1}=b_1.\end{subarray}}
N\big(\mathrm{x}^{E^{(1)}_i}=u_{1i}, 1\leq i\leq r_1\big). \ \ \ \ \ \ \ \ \ \ \ \ \ \  \ \ \ \ \ \ \ \ \ \ \ \ \ \ \ \ (3.12)
\end{align*}

Now we treat with the sum
$$
\sum_{\begin{subarray}{I}(u_{11},...,u_{1r_1})\in (\mathbb{ F}_{q}^*)^{r_1},\\
a_{11}u_{11}+...+a_{{1r_1}}u_{1r_1}=b_1.\end{subarray}}
N\big(\mathrm{x}^{E^{(1)}_i}=u_{1i}, 1\leq i\leq r_1\big).
\ \ \ \ \ \ \ \ \ \ \ \ \ \ \ \ \ \ \ \ \ \ \ \ \ \ \ \ \ \ \ \ \  \ \ \ \ \ \ \ \ \ \ (3.13)
$$
For any given $(u_{11},...,u_{1r_1})\in (\mathbb{ F}_{q}^*)^{r_1}$
with $\sum_{i=1}^{r_1}a_{1i}u_{1i}=b_1$, one has that
$$N\big(\mathrm{x}^{E^{(1)}_i}=u_{1i}, 1\leq i\leq r_1\big)$$
equals the number of the solutions
$(x_1,..., x_{n_1})\in (\mathbb{F}^*_q)^{n_1}$
of the following system of equations:
$$\left\{
\begin{aligned}
x_1^{e^{(1)}_{11}}...x_{n_1}^{e^{(1)}_{1n_1}}&=u_{11}, \\
............ \\
x_1^{e^{(1)}_{r_1,1}}...x_{n_1}^{e^{(1)}_{r_1,n_1}}&=u_{1r_1}.
\end{aligned}
\right. \eqno(3.14)
$$
Since $u_{1i}\neq0$ $(1\leq i\leq r_1)$, we can  deduce that the number of the solutions
$(x_1,..., x_{n_1})\in (\mathbb{F}^*_q)^{n_1}$ of (3.14) is equal to
the number of nonnegative integral solutions
$({\rm ind}_{\alpha}(x_1), ...,{\rm ind}_{\alpha}(x_{n_1}))\in \mathbb{N}^{n_1}$
of the following system of congruences
$$\left\{
\begin{aligned}
\sum\limits_{i=1}^{n_1}e^{(1)}_{1i}\text{ind}_\alpha(x_i)
&\equiv\text{ind}_\alpha (u_{11})\pmod{q-1}, \\
......&...... \\
\sum\limits_{i=1}^{n_1}e^{(1)}_{r_1,i}\text{ind}_\alpha(x_i)
&\equiv\text{ind}_\alpha (u_{1r_1})\pmod{q-1}.
\end{aligned}
\right. \eqno(3.15)
$$
But Lemma 2.2 tells us that (3.15) has
solutions $({\rm ind}_{\alpha}(x_1), ..., {\rm ind}_{\alpha}(x_{n_1}))
\in \mathbb{N}^{n_1}$ if and only if the extra conditions $(\widetilde{SE}^{(n_1)})$ hold.
Further, Lemma 2.2 gives us the number of
solutions $({\rm ind}_{\alpha}(x_1), ..., {\rm ind}_{\alpha}(x_{n_1}))
\in \mathbb{N}^{n_1}$ of (3.15) which is equal to
$$(q-1)^{n_1-s_1}\prod\limits_{i=1}^{s_1}\text{gcd}(q-1,d_i^{(1)}).$$
Hence
$$
N\big(\mathrm{x}^{E^{(1)}_i}=u_{1i}, 1\leq i\leq r_1\big)=(q-1)^{n_1-s_1}
\prod\limits_{i=1}^{s_1}\gcd(q-1, d_i^{(1)}).\eqno(3.16)
$$
Since
$$
H_1=
\sum_{\begin{subarray}{I}(u_{11},...,u_{1r_1})\in (\mathbb{F}_{q}^*)^{r_1},\\
a_{11}u_{11}+...+a_{1r_1}u_{1r_1}=b_1\\
 \ {\rm and} \ (\widetilde{SE}^{(n_1)}) \ {\rm holds.}\end{subarray}}1,
$$
it then follows from (3.16) that
\begin{align*}
&\sum_{\begin{subarray}{I}(u_{11},...,u_{1r_1})\in (\mathbb{ F}_{q}^*)^{r_1},\\
a_{11}u_{11}+...+a_{{1r_1}}u_{1r_1}=b_1.\end{subarray}}
N\big(\mathrm{x}^{E^{(1)}_i}=u_{1i}, 1\leq i\leq r_1\big) \\
&=\sum_{\begin{subarray}{I}(u_{11},...,u_{1r_1})\in (\mathbb{ F}_{q}^*)^{r_1},\\
a_{11}u_{11}+...+a_{1r_1}u_{1r_1}=b_1\\
 \ {\rm and} \ (\widetilde{SE}^{(n_1)}) \ {\rm holds.}\end{subarray}}
(q-1)^{n_1-s_1}\prod\limits_{i=1}^{s_1}\gcd(q-1, d_i^{(1)})\\
&=(q-1)^{n_1-s_1}\prod\limits_{i=1}^{s_1}\gcd(q-1,d_i^{(1)})
\sum_{\begin{subarray}{I}(u_{11},...,u_{1r_1})\in (\mathbb{ F}_{q}^*)^{r_1},\\
a_{11}u_{11}+...+a_{1r_1}u_{1r_1}=b_1\\
 \ {\rm and} \ (\widetilde{SE}^{(n_1)}) \ {\rm holds.}\end{subarray}}1\\
&=H_1(q-1)^{n_1-s_1}\prod\limits_{i=1}^{s_1}\gcd(q-1,d_i^{(1)}). \ \ \ \ \ \ \ \ \ \ \ \  \ \ \ \ \ \ \ \
\ \ \ \ \ \ \ \ \ \ \ \  \ \ \ \ \ \ \ \ \ \ \ \ \ \ \ \ \ (3.17)
\end{align*}
Thus by (3.6), (3.12) and (3.17), we obtain that
$M_{1}=N_1$ if $n_1<n_3<n_2$ and $b_2=0$.

By the same argument, one can deduce that either  $n_3>n_2$ or $n_1=n_3$  we have $M_{1}=N_1$.
This ends the second part of the Lemma 3.1.

Now we treat the third part of the Lemma 3.1. From (3.3) and the definition of $T(V_2)$,
we know that for any vector $v\in T(V_2)$ if and only if $n_3<n_2$.
Thus by the definition of $N_2$ and (3.4), we get that $M_2=N_2=0$ if $n_3\geq n_2.$

Let $n_3<n_2$. If $n_1<n_3<n_2<n_4$, it then follows from (3.4) and the definition of $T(V_2)$ that
\begin{align*}
&M_2=\sum_{\begin{subarray}{I}(u_{11},...,u_{1r_1},u_{21},...,u_{2r_3})\in (\mathbb{F}_q^*)^{r_1+r_3},\\
a_{11}u_{11}+...+a_{1r_1}u_{1r_1}=b_1,\\
a_{21}u_{21}+...+a_{2r_3}u_{2r_3}=b_2.
\end{subarray}}
N\Big(\left\{
\begin{aligned} \mathrm{x}^{E^{(1)}_i}&=u_{1i},\ 1\leq i\leq r_1,\\
\mathrm{x}^{E^{(1)}_j}&=0,\ r_1+1\leq j\leq r_2,\\
\mathrm{x}^{E^{(2)}_i}&=u_{2i},\ 1\leq i\leq r_3,\\
\mathrm{x}^{E^{(2)}_k}&=0,\ r_3+1\leq k\leq r_4.
\end{aligned}
\right.\Big).
 \ \ (3.18)
\end{align*}
The definitions of $\mathrm{x}^{E^{(1)}_j}$ and 
$\mathrm{x}^{E^{(2)}_k}$ tell us that
the fact that $\mathrm{x}^{E^{(1)}_j}=0$ and $\mathrm{x}^{E^{(2)}_k}=0$ 
with $r_1+1\le j\le r_2$, $r_3+1\leq k\leq r_4$
is reduced to saying that $\mathrm{x}^{E^{(1)}_{r_1+1}}=0$. 
It then follows from (3.18) that
\begin{align*}
M_2=\sum_{\begin{subarray}{I}(u_{11},...,u_{1r_1},u_{21},...,u_{2r_3})\in (\mathbb{F}_q^*)^{r_1+r_3},\\
a_{11}u_{11}+...+a_{1r_1}u_{1r_1}=b_1,\\
a_{21}u_{21}+...+a_{2r_3}u_{2r_3}=b_2.
\end{subarray}}
N\Big(\left\{
\begin{aligned} &\mathrm{x}^{E^{(1)}_i}=u_{1i},\ 1\leq i\leq r_1,\\
&\mathrm{x}^{E^{(2)}_j}=u_{2j},\ 1\leq j\leq r_3,\\
&\mathrm{x}^{E^{(1)}_{r_1+1}}=0.
\end{aligned}
\right.\Big).
 \ \ (3.19)
\end{align*}

It is easy to see that $\mathrm{x}^{E^{(1)}_{r_1+1}}=0$ is equivalent to $x_1...x_{n_2}=0$.
Since $n_3<n_2$, $u_{2r_3}\ne 0$ and $u_{2r_3}= x_1^{e^{(2)}_{r_3,1}}
...x_{n_3}^{e^{(2)}_{r_3,n_3}}$, one has $x_1...x_{n_3}\ne 0$.
So $\mathrm{x}^{E^{(1)}_{r_1+1}}=0$ is equivalent to $x_{n_3+1}...x_{n_2}=0$.
Then by (3.19), one gets that
\begin{align*}
M_2=\sum_{\begin{subarray}{I}(u_{11},...,u_{1r_1},u_{21},...,u_{2r_3})\in (\mathbb{F}_q^*)^{r_1+r_3},\\
a_{11}u_{11}+...+a_{1r_1}u_{1r_1}=b_1,\\
a_{21}u_{21}+...+a_{2r_3}u_{2r_3}=b_2.
\end{subarray}}
N\Big(\left\{
\begin{aligned} &\mathrm{x}^{E^{(1)}_i}=u_{1i},\ 1\leq i\leq r_1,\\
&\mathrm{x}^{E^{(2)}_j}=u_{2j},\ 1\leq j\leq r_3,\\
&x_{n_3+1}...x_{n_2}=0.
\end{aligned}
\right.\Big).
 \ \ \ \ \ \ \ \ \ \ \ \ \ \ \ \ \ (3.20)
\end{align*}

For any given $(u_{11},...,u_{1r_1},u_{21},...,u_{2r_3})\in (\mathbb{ F}_{q}^*)^{r_1+r_3}$
with $\sum_{i=1}^{r_1}a_{1i}u_{1i}=b_1$ and $\sum_{j=1}^{r_3}a_{2j}u_{2j}=b_2$, one has
\begin{align*}
&N\Big(\left\{
\begin{aligned} &\mathrm{x}^{E^{(1)}_i}=u_{1i},\ 1\leq i\leq r_1,\\
&\mathrm{x}^{E^{(2)}_j}=u_{2j},\ 1\leq j\leq r_3,\\
&x_{n_3+1}...x_{n_2}=0.
\end{aligned}
\right.\Big)=\#\{(x_1,...,x_{n_4})\in (\mathbb{F}_q)^{n_4}: \left\{
\begin{aligned} &\mathrm{x}^{E^{(1)}_i}=u_{1i},\ 1\leq i\leq r_1,\\
&\mathrm{x}^{E^{(2)}_j}=u_{2j},\ 1\leq j\leq r_3,\\
&x_{n_3+1}...x_{n_2}=0.
\end{aligned}
\right.\}.
\end{align*}
Since each of the components $x_{n_2+1},..., x_{n_4}$ can run over
the whole finite field $\mathbb{F}_q$ independently, it then follows that
\begin{align*}
&N\Big(\left\{
\begin{aligned} &\mathrm{x}^{E^{(1)}_i}=u_{1i},\ 1\leq i\leq r_1,\\
&\mathrm{x}^{E^{(2)}_j}=u_{2j},\ 1\leq j\leq r_3,\\
&x_{n_3+1}...x_{n_2}=0.
\end{aligned}
\right.\Big)\\
&=q^{n_4-n_2}\times\#\{(x_1,...,x_{n_2})\in (\mathbb{F}_q)^{n_2}: \left\{
\begin{aligned} &\mathrm{x}^{E^{(1)}_i}=u_{1i},\ 1\leq i\leq r_1,\\
&\mathrm{x}^{E^{(2)}_j}=u_{2j},\ 1\leq j\leq r_3,\\
&x_{n_3+1}...x_{n_2}=0.
\end{aligned}
\right.\}.\ \ \ \ \ \ \ \ \ \ \ \ \ \ \ (3.21)
\end{align*}
Notice that the choice of $(x_1,...,x_{n_3})\in (\mathbb{F}_q^*)^{n_{3}}$
satisfying that $\mathrm{x}^{E^{(1)}_i}=u_{1i}$ ($i=1,...,r_1$) 
and $\mathrm{x}^{E^{(2)}_j}=u_{2j}$ ($j=1,...,r_3$)
is independent of the choice of
$(x_{n_3+1},...,x_{n_2})\in (\mathbb{F}_q)^{n_2-n_3}$
satisfying that $x_{n_3+1}...x_{n_2}=0$. We then derive that
\begin{align*}
&\#\{(x_1,...,x_{n_2})\in (\mathbb{F}_q)^{n_2}: \left\{
\begin{aligned} &\mathrm{x}^{E^{(1)}_i}=u_{1i},\ 1\leq i\leq r_1,\\
&\mathrm{x}^{E^{(2)}_j}=u_{2j},\ 1\leq j\leq r_3,\\
&x_{n_3+1}...x_{n_2}=0.
\end{aligned}
\right.\}\\
&=\#\{(x_1,...,x_{n_3})\in (\mathbb{F}_q)^{n_3}: \left\{
\begin{aligned} &\mathrm{x}^{E^{(1)}_i}=u_{1i},\ 1\leq i\leq r_1,\\
&\mathrm{x}^{E^{(2)}_j}=u_{2j},\ 1\leq j\leq r_3.
\end{aligned}
\right.\}\times \\
&\#\{(x_{n_3+1},...,x_{n_2})\in (\mathbb{F}_q)^{n_2-n_3}: x_{n_3+1}...x_{n_2}=0.\}.
\ \ \ \ \ \ \ \ \ \ \ \ \ \ \ \ \ \ \ \ \  (3.22)
\end{align*}

On the other hand, we can easily compute that
$$
\#\{(x_{n_3+1},...,x_{n_2})\in (\mathbb{F}_q)^{n_2-n_3}: x_{n_3+1}...x_{n_2}=0\}=q^{n_2-n_3}-(q-1)^{n_2-n_3}.
\eqno(3.23)
$$
So by (3.21) to (3.23), one obtains that
\begin{align*}
&N\Big(\left\{
\begin{aligned} &\mathrm{x}^{E^{(1)}_i}=u_{1i},\ 1\leq i\leq r_1,\\
&\mathrm{x}^{E^{(2)}_j}=u_{2j},\ 1\leq j\leq r_3,\\
&x_{n_3+1}...x_{n_2}=0.
\end{aligned}
\right.\Big)\\
&=q^{n_4-n_2}(q^{n_2-n_3}-(q-1)^{n_2-n_3})\times\\
&\#\{(x_1,...,x_{n_3})\in (\mathbb{F}_q)^{n_3}: \left\{
\begin{aligned} &\mathrm{x}^{E^{(1)}_i}=u_{1i},\ 1\leq i\leq r_1,\\
&\mathrm{x}^{E^{(2)}_j}=u_{2j},\ 1\leq j\leq r_3.
\end{aligned}
\right.\}\\
&=q^{n_4-n_2}(q^{n_2-n_3}-(q-1)^{n_2-n_3})N\Big(\left\{
\begin{aligned} &\mathrm{x}^{E^{(1)}_i}=u_{1i},\ 1\leq i\leq r_1,\\
&\mathrm{x}^{E^{(2)}_j}=u_{2j},\ 1\leq j\leq r_3.
\end{aligned}
\right.\Big).\ \ \ \ \ \ \ \ \ \ \ \ \ \ \ (3.24)
\end{align*}
Then by (3.20) together with (3.24), we have
\begin{align*}
&M_2=q^{n_4-n_2}(q^{n_2-n_3}-(q-1)^{n_2-n_3})\times \\
&\sum_{\begin{subarray}{I}(u_{11},...,u_{1r_1},u_{21},...,u_{2r_3})\in (\mathbb{F}_q^*)^{r_1+r_3},\\
a_{11}u_{11}+...+a_{1r_1}u_{1r_1}=b_1,\\
a_{21}u_{21}+...+a_{2r_3}u_{2r_3}=b_2.
\end{subarray}}
N\Big(\left\{
\begin{aligned} &\mathrm{x}^{E^{(1)}_i}=u_{1i},\ 1\leq i\leq r_1,\\
&\mathrm{x}^{E^{(2)}_j}=u_{2j},\ 1\leq j\leq r_3.
\end{aligned}
\right.\Big). \ \ \ \ \ \ \ \ \ \ \ \ \ \  \ \ \ \  (3.25)
\end{align*}

Now the argument for (3.13) applied to the sum
$$
\sum_{\begin{subarray}{I}(u_{11},...,u_{1r_1},u_{21},...,u_{2r_3})\in (\mathbb{F}_q^*)^{r_1+r_3},\\
a_{11}u_{11}+...+a_{1r_1}u_{1r_1}=b_1,\\
a_{21}u_{21}+...+a_{2r_3}u_{2r_3}=b_2.
\end{subarray}}
N\Big(\left\{
\begin{aligned} &\mathrm{x}^{E^{(1)}_i}=u_{1i},\ 1\leq i\leq r_1,\\
&\mathrm{x}^{E^{(2)}_j}=u_{2j},\ 1\leq j\leq r_3.
\end{aligned}
\right.\Big)
$$
gives us that

$$
M_{2}=q^{n_4-n_2}(q^{n_2-n_3}-(q-1)^{n_2-n_3})
H_3(q-1)^{n_3-s_3}\prod\limits_{i=1}^{s_3}\gcd(q-1,d_i^{(3)})=N_2.
$$

This concludes the proof of the third part of the Lemma 3.1.

Next we prove the fourth part of Lemma 3.1.
The definition of $T(V_3)$ together with (3.3) give us
that for any vector $v\in T(V_3)$ if and only if
 $n_1\leq n_3<n_2<n_4$ or $n_2\leq n_3$.
Using the same argument of the third part of Lemma 3.1,
one can get that either $n_1\leq n_3<n_2<n_4$ or 
$n_2\leq n_3$, $M_3=N_3$. For the other case, one has $M_3=N_3=0$.

Now we deal with the fifth part of Lemma 3.1. 
It follows (3.4) and the definition of $T(V_4)$
that
$$
M_4=\sum_{\begin{subarray}{I}(u_{11},...,u_{1r_2},
u_{21},...,u_{2r_4})\in (\mathbb{F}_q^*)^{r_2+r_4},\\
a_{11}u_{11}+...+a_{1r_1}u_{1r_2}=b_1,\\
a_{21}u_{21}+...+a_{2r_3}u_{2r_4}=b_2.
\end{subarray}}
N\Big(\left\{
\begin{aligned} &\mathrm{x}^{E^{(1)}_i}=u_{1i},\ 1\leq i\leq r_2,\\
&\mathrm{x}^{E^{(2)}_j}=u_{2i},\ 1\leq i\leq r_4.
\end{aligned}
\right.\Big).  \ \ \ \ \ \ \ \ \ \ \ (3.26)
$$

If $n_2\leq n_4$, the argument of (3.13) applied to (3.26) tells us that
$$M_4=H_4(q-1)^{n_4-s_4}\prod\limits_{j=1}^{s_4}\gcd(q-1,d_j^{(4)})=N_4.$$

Similarly, if $n_2> n_4$, we can get that
$$M_4=H_2(q-1)^{n_2-s_2}\prod\limits_{j=1}^{s_2}\gcd(q-1,d_j^{(2)})=N_4.$$

Next we prove the sixth part of Lemma 3.1. The definition of $T(V_5)$ and (3.3) give us that
the vector $v\in T(V_5)$ if and only if $n_2>n_4$. If $n_2\leq n_4$,
by the definition of $M_5$, one can easily deduce that $M_5=N_5=0$. Now we let $n_2>n_4$.
It follows from (3.4) and the definition of $T(V_5)$ that
\begin{align*}
M_5=\sum_{\begin{subarray}{I}(u_{11},...,u_{1r_1},u_{21},...,u_{2r_4})\in (\mathbb{F}_q^*)^{r_1+r_4},\\
a_{11}u_{11}+...+a_{1r_1}u_{1r_1}=b_1,\\
a_{21}u_{21}+...+a_{2r_4}u_{2r_4}=b_2.
\end{subarray}}
N\Big(\left\{
\begin{aligned} \mathrm{x}^{E^{(1)}_i}&=u_{1i},\ 1\leq i\leq r_1,\\
\mathrm{x}^{E^{(1)}_j}&=0,\ r_1+1\leq j\leq r_2,\\
\mathrm{x}^{E^{(2)}_k}&=u_{2k},\ 1\leq k\leq r_4.
\end{aligned}
\right.\Big).
 \ \ (3.27)
\end{align*}
Using the same argument of (3.19), we obtain that
$$ M_5=\big(q^{n_2-n_4}-(q-1)^{n_2-n_4}\big)H_4(q-1)^{n_4-s_4}\prod\limits_{j=1}^{s_4}\gcd(q-1,d_j^{(4)})=N_5.$$

In the rest, we prove $M_6=N_6$.
From the definition of $T(V_6)$ and (3.3), we konw that
the vector $v\in T(V_6)$ if and only if $n_3>n_2$. Thus one has $M_6=N_6=0$, if $n_3\leq n_2$.
Now we let $n_3>n_2$. Then (3.4) and the definition of $T(V_6)$ tell us that
\begin{align*}
M_6
=\sum_{\begin{subarray}{I}(u_{11},...,u_{1r_2})\in (\mathbb{F}_q^*)^{r_2},\\
a_{11}u_{11}+...+a_{1r_2}u_{1r_2}=b_1\ {\rm {and} }
\ b_2=0.
\end{subarray}}
N\Big(\left\{
\begin{aligned}  \mathrm{x}^{E^{(1)}_i}&=u_{1i},\ 1\leq i\leq r_2,\\
 \mathrm{x}^{E^{(2)}_j}&=0,\ 1\leq j\leq r_4.
\end{aligned}
\right.\Big).
 \ \ (3.28)
\end{align*}
Using the same argument of (3.5), we get that
\begin{align*}
M_6={\left\{\begin{array}{rl}
q^{n_4-n_3}(q^{n_3-n_2}-(q-1)^{n_3-n_2})
H_2(q-1)^{n_2-s_2}&\prod\limits_{j=1}^{s_2}\gcd(q-1,d_j^{(2)}),\\
& {\rm if} \ n_3>n_2 \ {\rm and} \ b_2=0,\\
0,   & {\rm otherwise},
\end{array}\right.}
\end{align*}
as desired. This completes the proof of the Lemma 3.1. \hfill$\Box$

We can now use Lemma 3.1 to prove Theorem 1.3 as the conclusion of this section.

{\bf Proof of Theorem 1.3.} We divide the proof into the following eight cases.

{\it Case 1.} $n_1<n_3<n_2<n_4$.
From (3.3) we can obtain that
if the vectors $v\in T$ have nonzero components of the vectors
$(u_{11},...,u_{1r_2},u_{21},...,u_{2r_4})\in T$, then
we must have $v\in v_1$, or $v\in v_2$, or $v\in v_3$, or $v\in v_4$.
Then we have $T=\bigcup\limits_{i=0}^{4}T(V_i)$.
It follows from (3.2) and (3.4) that $N=\sum\limits_{i=0}^4M_i$.
So by the Lemma 3.1, we have $N=\sum\limits_{i=0}^4N_i$ as desired.

{\it Case 2.} $n_2=n_3$.
Since (3.3) tells us that if the vectors $v\in T$
have nonzero components of the vectors
$(u_{11},...,u_{1r_2},u_{21},...,u_{2r_4})\in T$, then
we must have $v\in v_1$, or $v\in v_3$, or $v\in v_4$.
Then $T=\bigcup\limits_{i=0\atop i\neq2}^{4}T(v_i).$
Thus by (3.2), (3.4) and the Lemma 3.1, the desired result
$N=\sum\limits_{i=0\atop i\neq2}^4N_i$ follows.

{\it Case 3.} $n_1=n_3$, $n_2<n_4$.
Similarly (3.3) tells us that if the vectors $v\in T$ 
have nonzero components of the vectors
$(u_{11},...,u_{1r_2},u_{21},...,u_{2r_4})\in T$, then
we must have $v\in v_2$, or $v\in v_3$, or $v\in v_4$.
Thus $T=\bigcup\limits_{i=0\atop i\neq1}^{4}T(v_i)$.
It follows from (3.2), (3.4) and the Lemma 3.1 that 
$N=\sum\limits_{i=0\atop i\neq1}^4N_i.$

{\it Case 4.} $n_1<n_3$, $n_2=n_4$.
It follows from (3.3) and the vectors $v\in T$ have 
nonzero components of the vectors
$(u_{11},...,u_{1r_2},u_{21},...,u_{2r_4})\in T$ that
$v\in v_1$, or $v\in v_2$, or $v\in v_4$.
So $T=\bigcup\limits_{i=0\atop i\neq3}^{4}T(v_i).$
Then by (3.2), (3.4) and the Lemma 3.1, one gets that
$N=\sum\limits_{i=0\atop i\neq3}^4N_i$.

{\it Case 5.} $n_1=n_3$, $n_2=n_4$. By the same argument,
we know that if the nonzero vectors $v\in T$, then $v\in v_2$, or $v\in v_4$.
It infers that $T=\bigcup\limits_{i=0\atop i\neq1,3}^{4}T(v_i)$.
Using (3.2), (3.4) and the Lemma 3.1, we deduce that
$N=\sum\limits_{i=0\atop i\neq1, 3}^4N_i$.

{\it Case 6.} $n_1<n_3$, $n_2>n_4$. Similarly, from (3.3) we deduce that
if the vectors $v\in T$ have nonzero components of the vectors
$(u_{11},...,u_{1r_2},u_{21},...,u_{2r_4})\in T$, then
we obtain $v\in v_1$, or $v\in v_2$, or $v\in v_4$, or $v\in v_5$.
Thus $T=\bigcup\limits_{i=0\atop i\neq3}^{5}T(v_i)$.
From (3.2), (3.4) and the Lemma 3.1, the desired result
$N=\sum\limits_{i=0\atop i\neq3}^{5}N_i$ follows immediately.

{\it Case 7.} $n_1=n_3$, $n_2>n_4$. It's easy to deduce that
if the vectors $v\in T$ have nonzero components of the vectors
$(u_{11},...,u_{1r_2},u_{21},...,u_{2r_4})\in T$, then $v\in v_2$, or $v\in v_4$, or $v\in v_5$.
Then $T=\bigcup\limits_{i=0\atop i\neq1,3}^{5}T(v_i)$.
Using (3.2), (3.4) and the Lemma 3.1, we have $N=\sum\limits_{i=0\atop i\neq1,3}^{5}N_i$ as desired.

{\it Case 8.} $n_3>n_2$. From (3.3), one knows that if
the nonzero vectors $v\in T$, then  $v\in v_1$, or $v\in v_3$, or $v\in v_4$, or $v\in v_6$.
So $T=\bigcup\limits_{i=0\atop i\neq2,5}^{6}T(v_i)$.
Then by (3.2), (3.4) and the Lemma 3.1, the desired result $N=\sum\limits_{i=0\atop i\neq2,5}^{6}N_i$ follows.

Therefore Theorem 1.3 follows immediately. This ends the proof of Theorem 1.3. \hfill$\Box$

\section {An example}

In this section, we supply an example to illustrate
the validity of our main result.\\
{\bf Example 4.1.} We use Theorem 1.3 to compute the number $N$ of
rational points on the following variety over $\mathbb{F}_7$:
\begin{align*}
{\left\{\begin{array}{rl}
&x_1x_2^4+x_1x_2^5+x_1^2x_2^3x_3x_4^4-2=0,\\
&x_1x_2^5x_3^3+x_1x_2^3x_3^2+x_1^2x_2^4x_3^3x_4^5x_5x_6-4=0.\\
\end{array}\right.}
\end{align*}
Obviously, we have $n_1=2$, $n_2=4$, $n_3=3$, $n_4=6$,
$$E^{(1)}=\left( \begin{array}{*{20}c}
1&4\\
1&5\\
\end{array}\right),
E^{(2)}=\left( \begin{array}{*{20}c}
1&4&0&0\\
1&5&0&0\\
2&3&1&4\\
1&5&3&0\\
1&3&2&0\\
\end{array}\right),
$$

$$E^{(3)}=\left( \begin{array}{*{20}c}
1&4&0\\
1&5&0\\
1&5&3\\
1&3&2\\
\end{array}\right) \  \text{and} \
E^{(4)}=\left( \begin{array}{*{20}c}
1&4&0&0&0&0\\
1&5&0&0&0&0\\
2&3&1&4&0&0\\
1&5&3&0&0&0\\
1&3&2&0&0&0\\
2&4&3&5&1&1\\
\end{array}\right).$$
It then follows from the assumption $n_1<n_3<n_2<n_4$,
$2$, $4\in {F}_7^*$ and Theorem 1.3 that
$$
N=\sum\limits_{i=2}^{4}N_i.    \eqno(4.1)
$$

We first calculate $N_2$.
Using elementary transformations, One can easily deduce that
$$U^{(3)}=\left( \begin{array}{*{20}c}
1&0&0&0\\
-1&1&0&0\\
2&-2&1&-1\\
-6&5&-2&3
\end{array}\right)
\ {\rm and} \
V^{(3)}=\left( \begin{array}{*{20}c}
1&-4&0\\
0&1&0\\
0&0&1
\end{array}\right)$$
such that
$$U^{(3)}E^{(3)}V^{(3)}
={\rm SNF}(E^{(3)})=\left( \begin{array}{*{20}c}
1&0&0\\
 0&1&0\\
 0&0&1\\
 0&0&0
\end{array}\right).$$
Thus $d_1^{(3)}=d_2^{(3)}=d_3^{(3)}=1$ and $s_3=3$.
But Lemma 2.4 tells us that the number of the vectors
$(u_{11},u_{12},u_{21},u_{22})\in(\mathbb{F}^*_7)^4$ satisfying that
\begin{align*}
{\left\{\begin{array}{rl}
u_{11}+u_{12}=2\\
u_{21}+u_{22}=4
\end{array}\right.}
\end{align*}
is 25. For example, the vector
$$(u_{11},u_{12},u_{21},u_{22})=(3,6,1,3)$$
is one of the solutions.
Choose the primitive element 3 of $\mathbb{F}^*_7$. Then
\begin{align*}
(h'_1, h'_2,h'_3, h'_4)^T&=U^{(3)}(\text{ind}_33,\text{ind} _36, 
\text{ind} _31,\text{ind}_33)^T\\
&=U^{(3)}(1,3,6,1)^T\equiv(1,2,1,0)^T \pmod6.
\end{align*}
We deduce that the condition $(\widetilde{SE}^{(3)})$ 
that $1|1$, $1|2$ and $6|0$ holds.
Similarly, using Matlab we compute that $H_3=4$. So
\begin{align*}
N_2&=q^{\max\{n_2,n_4\}-\min\{n_2,n_4\}}(q^{\min\{n_2,n_4\}-n_3}-(q-1)^{\min\{n_2,n_4\}-n_3})L_3\\
&=q^2(q-(q-1))H_3(q-1)^{n_3-s_3}\prod\limits_{j=1}^{s_3}\gcd(q-1,d_j^{(3)})=196.
\end{align*}

Consequently, we turn our attention to the computation of $N_3$.
Using the elementary transformations, one gets that
$$U^{(2)}=\left( \begin{array}{*{20}c}
1&0&0&0&0\\
-1&1&0&0&0\\
-7&5&1&0&0\\
-9&7&1&-1&1\\
-6&5&0&-2&3
\end{array}\right)
\ {\rm and} \
V^{(2)}=\left( \begin{array}{*{20}c}
1&-4&0&0\\
0&1&0&0\\
0&0&1&-4\\
0&0&0&1
\end{array}\right)$$
such that
$$U^{(2)}E^{(2)}V^{(2)}
={\rm SNF}(E^{(2)})=\left( \begin{array}{*{20}c}
1&0&0&0\\
0&1&0&0\\
0&0&1&0\\
0&0&0&4\\
0&0&0&0
\end{array}\right).$$
Thus $d_1^{(2)}=d_2^{(2)}=d_3^{(2)}=1$, $d_4^{(2)}=4$ and $s_2=4$.
One can easily deduce that the number of the vectors
$(u_{11},u_{12},u_{13},u_{21},u_{22})\in(\mathbb{F}^*_7)^5$
such that
\begin{align*}
{\left\{\begin{array}{rl}
u_{11}+u_{12}+u_{13}=2\\
u_{21}+u_{22}=4
\end{array}\right.}
\end{align*}
is 155. Choose the primitive element 3 of $\mathbb{F}^*_7$.
The argument for calculating $H_3$ and using Matlab
we compute that $H_2=9$. Hence
\begin{align*}
N_3&=\big(q^{n_4-n_2}-(q-1)^{n_4-n_2}\big)L_2 \\
&=\big(q^{n_4-n_2}-(q-1)^{n_4-n_2}\big)H_2(q-1)^{n_2-s_2}\prod\limits_{j=1}^{s_2}
\gcd(q-1,d_j^{(2)})\\
&=18\big(q^2-(q-1)^2\big)=234.
\end{align*}

Let us now calculate $N_4$. Using the elementary transformations, we obtain that
$$U^{(4)}=\left( \begin{array}{*{20}c}
1&0&0&0&0&0\\
-1&1&0&0&0&0\\
-7&5&1&0&0&0\\
15&-11&-3&0&0&1\\
-9&7&1&-1&1&0\\
-6&5&0&-2&3&0
\end{array}\right)$$
and
$$V^{(4)}=\left( \begin{array}{*{20}c}
1&-4&0&0&0&0\\
0&1&0&0&0&0\\
0&0&1&0&-4&0\\
0&0&0&0&1&0\\
0&0&0&1&7&-1\\
0&0&0&0&0&1
\end{array}\right)$$
such that
$$U^{(4)}E^{(4)}V^{(4)}
={\rm SNF}(E^{(4)})=\left( \begin{array}{*{20}c}
1&0&0&0&0&0\\
0&1&0&0&0&0\\
0&0&1&0&0&0\\
0&0&0&1&0&0\\
0&0&0&0&4&0\\
0&0&0&0&0&0
\end{array}\right).$$
Thus $d_1^{(4)}=d_2^{(4)}=d_3^{(4)}=d_4^{(4)}=1$,
$d_5^{(4)}=4$ and $s_4=5$. So Lemma 2.4 gives us 
that the number of the vectors
$(u_{11},u_{12},u_{13},u_{21},u_{22},u_{23})
\in (\mathbb{F}^*_7)^6$ such that
\begin{align*}
{\left\{\begin{array}{rl}
u_{11}+u_{12}+u_{13}=2\\
u_{21}+u_{22}+u_{23}=4
\end{array}\right.}
\end{align*}
is equal to 961. Choose the primitive element 3 of
$\mathbb{F}^*_7$. By the argument for calculating $H_3$
and using Matlab, we compute that $H_4=84$. Thus one has
\begin{align*}
N_4&=l_4=H_4(q-1)^{n_4-s_4}\prod\limits_{j=1}^{s_4}\gcd(q-1,d_j^{(4)})\\
&=84\times 6 \times 2=1008.
\end{align*}
Finally, by (4.1), we have
$$N=\sum\limits_{i=2}^{4}N_i=196+234+1008=1438.$$

%%%%%%%%%%%%%%%%%%%%%%%%%%%%%%%%%%%%%%%%%%%%%%%%%%%%%%%%%%%%%%%%%%%%%%%%%%%%%%%%%%%%%%

\end{document}